\documentclass[12pt]{amsart}
\newtheorem{thm}{Theorem}[section]
\newtheorem{cor}[thm]{Corollary}
\newtheorem{lem}[thm]{Lemma}
\newtheorem{prop}[thm]{Proposition}
\theoremstyle{definition}
\newtheorem{defn}[thm]{Definition}
\theoremstyle{remark}
\newtheorem{rem}[thm]{Remark}
\numberwithin{equation}{section}
\newcommand{\Real}{\mathbb R}

\newcommand{\eps}{\varepsilon}

\begin{document}

\title[]{Topological equivalence  of smooth functions with
isolated critical points on a closed surface}
\author{A.O.Prishlyak}
\address{Kiev University,  Department of Geometry}
\email{prish@mechmat.univ.kiev.ua, topology@imath.kiev.ua}



\begin{abstract}
We consider functions with isolated critical points on a closed
surface. We prove that in a neighborhood of a critical point
the function conjugates with Re$z^k$ for the some nonnegative
integer $k$. The full topological invariant
of such functions is constructed.
\end{abstract}

\maketitle

Let $F$ be a closed smooth surface, $f, g: F \to \Real$ be a
smooth functions with a finite number of critical points. In view
of the compactness of the surface, this condition is equivalent
that each critical point is isolated.

Functions $f, g$ are called topologically equivalent if
there are homeomorphisms $h: F \to F$ and $h': \Real \to \Real$
such that $f h = h' g.$ We say that the functions are
topologically conjugate if they are topologically equivalent
and homeomorphism $h'$ preserves the orientation. Homeomorphisms
h and $h'$ we call conjugate homeomorphisms.

In the theory of singularities the differentiable equivalence
are studied. It is such a topological equivalence, in which the
conjugate homeomorphisms are diffeomorphisms.

The purpose of the given paper is to give a topological
classification of smooth functions on a closed surface with
isolated critical points. For Morse function such classification
was obtained by Sharko [5] and Kulinich [1]. There is a
classification for functions with unique critical level (except
minimum and maximum) that is imbedded graph [2].

In the section 1 it is proved that in a neighborhood of an isolated
critical point the function is topologically
conjugated with the function Re$z^k$, where $k$ is a
nonnegative integer. In the section 2 the full
topological invariant of the function on a closed surface is
constructed and the criterions of topological conjugation of
functions with isolated singular points are proved.
In the section 3 the numbers of minimal non-conjugated
functions and the numbers of minimal topologically
non-equivalent functions are calculated for closed oriented
and non-oriented surfaces of genus 0, 1, 2 and 3.

\section{Local conjugation of functions in a neighborhood of
an isolated critical point.}

Let $F$ be a smooth closed surface and $f:F\to \Real$ be a
smooth function.

\begin{lem}
For each isolated critical point $x_0$ with
the critical value $y_0=f(x_0)$ there exists a neighborhood
$U(x_0),$ in which

$$f^{-1}(y_0) \cap (U (x_0)) \simeq  Con (\cup S^0).$$

\end{lem}

Here $Con (\cup S^0)$ is a cone on the union of 0-dimensional
spheres $S^0$, that
is the union of the even number of the segments with a unique common
point, which is an endpoint for each of the segments.

\begin{proof}
 Since $x_0$ is an isolated singular point, there is
a simply connected neighborhood $U,$ which closure does not
contain other critical points of the function and which
boundary is a smooth closed curve. Let us assume $g$ is a
restriction of function $f$ to the set $U\setminus \{x_0\}.$ Then
the levels of both functions coincide in this set. As $y_0$ is
a regular value of the function g, the pre-image $g^{-1} (y_0)$
is a 1-dimensional submanifold in $U\setminus  \{x_0\},$ that
is a disjoint union of imbedded circles and open intervals.

Let $H$ be a  connected component  of the set $g^{-1}(y_0).$
If H is homeomorphic to a circle, then it is the boundary of
a closed set $D$. As function $f$ is not constant on $D$ the
maximum or minimum value of it on $D$ is distinct from $y_0,$
and the appropriate point of maximum or minimum  is another,
except for $x_0,$ critical point in domain $U.$ The obtained
contradiction proves that $H$ is homeomorphic to an open
interval.

Let us show that all accumulation points of $H$ in
$U\setminus \{x_0\}$ belong to $H$. By contradiction, let
there is such a point $z$ that does not belong to $H$.
Then, according to the continuity of the function $f$:
$z \in g^{-1} (y_0).$ As $z\notin H$, then $z$ belong other
component $H_1$ of $g^{-1}(y_0).$ Since $g^{-1}(y_0)$ is
a submanifold, there is a neighborhood of a point $z$ that
does not contain other points from $g^{-1}(y_0),$
except the points of $H_1.$ It contradicts assumption  that $z$
is an accumulation point of $H.$

Using similar argument to the neighborhood of closure $cl(U)$,
we obtain that the accumulation points of $H$ in $cl(U)$
belong the same component of $f^{-1}(y_0) \cap cl(u)$ as $H$.
So in the boundary of $U$ there is no more than two accumulation
points of interval $H$. Accumulation points of $H$ that belong
to the boundary of $U$ or coincide with $x_0$ we call
extremities of the interval. Three types of the interval $H$ are
possible:

1) both extremities coincide with the point $x_0$,

2) both extremities lay on the boundary of $U$,

3) one of extremities coincides with a point $x_0$ and another
lays on the boundary U.

Let us prove that there doesn't exit a component of the level
$g^{-1}(y_0)$ that form a loop with vertex in the point $x_0$
(i.e. both extremities coincide with $x_0$).
By contradiction, let such component exists, then it is
a boundary of the two-dimensional disk. In this disk there
is a critical point of minimum or maximum that distinct from
$x_0.$ It contradicts a choice of a neighborhood $U(x_0).$ So
there doesn't exit an interval of type 1).

Note that  only the finite number of extremities of intervals
can lay on the boundary of $U$. Otherwise, the condition
that $f^{-1}(y_0)$ is a submanifold in a limit point is false.
So there are a finite number of intervals of type 2) and 3).

We contract a neighborhood $U$ to such a neighborhood
$U(x_0)$, which closure is not intersected with interval of
type 2) and such, that the boundary of $U(x_0)$ is smooth and
has transversal intersection with $ g^{-1}(y_0)$. Then
$f^{-1}(y_0) \cap cl(U (x_0))$
consists of finite union of the closed intervals, one of which
extremities coincides with a point $x_0,$ and another lays on
the boundary of $U(x_0)$. Hence $f^{-1}(y_0) \cap cl(U (x_0))$
is a cone of the finite number of points.

Let us prove that number of extremities laying on the boundary
$\partial U$ is even.
We call domains in $U(x_0) \setminus f^{-1}(y_0)$ with $f(x) >
y_0$ by positive and with $f(x)< y_0$ by negative. Each
interval contains in the boundary both positive and negative
domain. These domains alternate. Thus we have the same number of the
domains of each type. Hence, total
number of the domains, as well as intervals, is even.
\end{proof}

\begin{rem}
The similar lemma takes place for three-dimensional
manifolds [6].
\end{rem}

\begin{thm}
For each isolated critical point $x_0$
(except for the local minimum and maxima) of a smooth function
on a surface there is a neighborhood, in which the function
is topologically conjugated with the function Re$z^k$ for some
nonnegative integer $k.$
\end{thm}

\begin{proof}
Let $U(x_0)$ be a neighborhood with the same
properties as in the lemma and $2k$ be the number of arcs,
which an extremity is the point $x_0.$

The critical level devides the neighborhood $U(x_0)$ to domains.
Let $V$ be one of them and such that $f(x) < y_0$ for $x \in V.$
Then in the points of the  intersection
$\partial U(x_0)$ with $\partial V$ the vector field grad$f$
is directed inside of the domain $V.$

Let us prove that there is a trajectory that passes through
$\partial U(x_0) \cap \partial V$ and tends to $x_0.$
Really, let
$h: \partial U(x_0) \cap \partial V \to f^{-1}(y_0)$ be a
map given by the formula

$$h(x)=\begin{cases}
x_0 & \text{ if $\gamma (x)\cap f^{-1}(y_0)=\emptyset$},\\
\gamma (x)\cap f^{-1}(y_0) & \text{ if $\gamma (x)
\cap f^{-1}(y_0)\ne \emptyset$}.
\end{cases} $$

In view of the smoothness of the vector field grad$f$, the continuity of
map $h$ follows for points, which image does not coincide with
the point $x_0.$ If the set of such points coincides with
$\partial U(x_0) \cap \partial V$, then the image of this
connected set by the continuous map is connected set. The
obtained inconsistency proves the existence of a required
trajectory.

Let $W$ be a neighborhood  of a point $x_0$ in $f^{-1}(y_0)$.
We define a neighborhood $W(\eps)$ in F for $\eps >0$ by formula
$$ W(\eps) = \{x\in F: |f(x)-y_0| < \eps, \gamma (x) \cap W \ne
\emptyset\}.$$

By the same formula, as above we define a map
$$h: W(\eps) \to W.$$

Let us prove that for any connected neighborhood $W$ of the
point $x_0$ in $f^{-1}(y_0),$ such that $cl(W) \cap U(x_0)
\ne \emptyset$
there exists $\eps >0$ that $W(\eps) \subset U(x_0).$
Let $z_i \in \partial W, i=1, 2, ..., 2k.$ Let $z_i',
z_i'' \in \gamma (z_i) \cap \partial U(x_0)$ be such points
that $f(z_i')<y_0<f(z_i'')$ and between points $z_i',
z_i''$  on the trajectory $\gamma (z_i)$
there don't exist points of its intersection with boundary
$\partial U(x_0).$ Let $\eps _i=min \{|y_0-f(z_i')|, |y_0-
f(z_i")|\}.$ Denote by $T$ union of all arcs of the boundary
$\partial U(x_0)$ that their endpoints belong to the set
$\{z_1', z_1'', z_2', z_2'', ..., z_{2k}', z_{2k}''\}$
and which do not contain  points of $f^{-1}(y_0)$. Let
$\eps < min_{z\in T} \{\eps_i, |y_0-f(z)|\}.$ Then it is easy
to see that $W(\eps) \subset  U(x_0).$

Let us construct an conjugated  homeomorphism from this
neighborhood to appropriate neighborhood of the function
Re$z^k.$ We  take a set
$$\{ z: |z | < 1, arg z =  (2n+1) /2k, 0\leq n \leq k-1,
n\in N\}$$
in the capacity of the neighborhood W of the function Re$z^k.$
Let us select one-to-one the correspondence of domains,
into which level lines decompose neighborhoods $W(\eps)$, so that
to adjacent domains corresponded adjacent ones and positive domains
to positive ones. We construct an conjugated homeomorphism
between appropriate domains. Let for a determination, domain $V$
corresponds domain with a polar angle:

$$ \pi / (2k) < \varphi  < 3\pi/(2k).$$

Let ${x_1, x_1'} = \partial V \cap \partial W,
x2 =\gamma (x_1) \cap f^{-1}(y_0-\eps),$ and $ x_3 $ be a
nearest to $x_2$ on $ \partial V \cap f^{-1}(y_0-\eps)$ point,
for which $h(x_3) =x_0.$ We fix a homeomorphism $p$ of the
curve  $\theta \subset  f^{-1} (y_0)$ with endpoints $x_0$ and
$x_1$ to $[0, 1]$ so, that $p (x_0) =0, p (x_1) =1.$ Let us
define map $H$ from the closed domain D, limited by a curvilinear
tetragon with vertex  $x_0, x_1, x_2, x_3$ and sides laying on
$f^{-1}(y_0), \gamma (x_2), f^{-1} (y_0-\eps), \gamma (x_3)$
to the rectangle $[0,1]\times [0,\eps]$ by formula
$$H(x) = (p (h (x)), f(x_0) -y_0).$$

Thus, the level lines of function and integrated curves are
mapped in segments parallel to the sides  of a rectangle.

Let us prove that the map $H$ is a homeomorphism.
Since each point lays on an integrated trajectory and the
function is monotonic at going on a trajectory in one direction,
the uniqueness of map $H$ follows. For an arbitrary point
$z\in D$ and sequence of points $\{z_i \in D\}$, it is
necessary to prove that $z_i \to z$ if and only if
$H(z_i) \to H(z)$.
By contradiction, let $z_i \to z$ and $H(z_i)$ does not
converge to $H(z).$ Then there is such $\eps _z>0$  and
subsubsequence $z_{i_k}$  that $|h(z_{i_k} )-h(z)|< \eps _z$
or $|f(z_{i_k})-f(z)| < \eps_z.$
It is equivalent to  that
$$ z_{i_k} \in  D \setminus (h^{-1}(h(z)-\eps _z, h(z)+ \eps _z)
\cap f^{-1}(f(z)-\eps_z, f(z)+\eps _z)).$$

Thus,   points $z_{i_k}$ lay out of the neighborhood of a point
$z,$ that are bounded by curves $h^{-1}(h(z)-\eps_z), \  f^{-1}
(f(z)-\eps_z), \ h^{-1}(h(z)+\eps_z), \  f^{-1}(f(z)+\eps_ z)).$
It contradicts with convergence of the sequence $z_i.$

Let us assume now that the sequence $H(z_i)$ converges to
$H(z),$ and the sequence $z_i$ does not converge to $z.$ Since
$D$ is the compact set then there is a subsubsequence $z_{i_k}$
that converges to $z_0 \ne z.$ Then $H(z_{i_k}) \to H(z_0).$ It
contradicts that $H(z_i) \to H(z).$

By analogy with above, we construct a curvilinear tetragon
with vertex in points $x_0, x_1', x_2', x_3'.$ If $x_3' =x_3,$
as well as above we create, a homeomorphism of this tetragon on
a rectangle $[0,1] \times [0,\eps].$ If $x_3' \ne x_3$ we denote
by $D'$ a curvilinear tetragon with vertex in the points
$x_0, x_1', x_2', x_3',$ and by $D"$ a curvilinear triangle with
vertex in points $x_0, x_3, x_3'.$ Let us fix homeomorphisms
$p, q$ of arcs on level lines of the function between the points
$x_0, x_1'$ and $x_3, x_3'$ to [0,1] and [0,1/2], accordingly.
Let us set a map $H:D'\cup D"\to [0,1]\times [0, \eps]$ by
the formula:

$$H(x)=\begin{cases}
(\ p(h(x))-{{1-p(h(x))} \over {2\varepsilon }}
f(x),\  -f(x)) & \text{ if $x\in D'$},\\

(\   -{{f(x)} \over \varepsilon } \ q(\gamma (x)\cap f^{-1}
(\varepsilon )),\  -f(x) ) & \text{ if $x\in D''$}.
\end{cases}
$$

As above, this map is a homeomorphism that maps a level lines
to horizontal curves.

Similarly we build homeomorphisms of domain for function Re$z^k.$
Then, having taken appropriate compositions of homeomorphisms,
we obtain homeomorphisms of domain of function $f$ to domain of
function Re$z^k.$ By the construction, these homeomorphisms
coincide in the boundaries of domain and thus set a required
homeomorphism of the neighborhood $W(\eps).$
\end{proof}

\begin{rem}
If $k=1,$ in a neighborhood
of the critical point, the level line of function is same as
well as in a neighborhood of a regular point (to within a
homeomorphism).
\end{rem}

\begin{rem}
Similarly to proof of the theorem one can show, that each local
minimum (maxima) has neighborhood, in which the function
conjugates with the functions $g(x, y) =x^2+y^2
(g (x, y) = -x^2-y^2).$
\end{rem}

\begin{defn}
 We call by a Poincare index  $ind_xf$ of a critical point
 $x$ of a function f a Poincare index of the gradient field
grad$f$ in some Riemannian metric.
\end{defn}

\begin{cor}
Let $x$ be an isolated critical point of a function $f,$ and
$y$ be a same of a function $g.$ The functions $f$ and $g$ are
topologically equivalent in some neighborhoods of these points
if and only if
$$ind_xf = ind_yg.$$
\end{cor}

\begin{rem}
For the function $f$ = Re$z^k$, the Poincare index   $ind_0f = 1-k.$
It is known [4] that for each critical point $x$ on a
surface $ind_xf\leq 1.$
\end{rem}

\section{Global conjugation of functions with isolated critical
points on a closed surface.}

\subsection {Diagram of the function.}

Let $F$ be a smooth closed surface and $f: F \to \Real$ be a smooth
function with isolated critical points and critical values
$y_1, y_2, ..., y_n$ that $y_i < y_j,$ if $i < j.$ From the
theorem 1 the pre-image $f^{-1} (y_i)$ of each critical value
$y_i$ is homeomorphic to the graph united with circles.
The edges of the graph and the circles are smoothly
imbedded except of critical points that coincide
with vertexes of the graph. On each circle we select one
point and we consider the circle as loop with a vertex in
the selected point.

Denote by $G_i(f)$  the graph (probably disconnected)
that coincide with a pre-image of a critical value $y_i$ of
the function $f.$ Vertexes of the graph are critical and selected
points.

\begin{defn}
The surface $F$ together with the graphs $G_i(f)$ that are
imbedded in it is called a diagram $D$ of the functions $f.$
So $D = \{F, G_1(f),..., G_n(f)\}$.  Two diagrams are called
isomorphic if there is a homeomorphism of surfaces, which maps
the graphs to the graphs, and the vertexes to then vertexes and
preserve the order of the graphs.
\end{defn}

\begin{thm}
Two functions with isolated critical points on closed
surfaces are topologically equivalent if and only if their
diagrams are isomorphic.
\end{thm}

\begin{proof}

Necessity. If a conjugated homeomorphism of the surface $F$ is
given, it maps critical levels to critical levels. It sets an
isomorphism of the diagrams.

Sufficiency. Suppose that the diagrams of two functions are
isomorphic. Then there is a homeomorphism $g$ of a surface $F$
that maps critical levels to critical levels. Let us cut
the surface $F$ by the critical levels. We obtain surfaces $F_k,$
which is homeomorphic to cylinders $S^1\times [0,1].$ Then the
homeomorphism $g$ induces homeomorphisms of obtained cylinders.
We replace these homeomorphisms with homeomorphisms, which
map the levels of the function to the levels proportionally to value of
function between two critical values and coincide with initial
ones on each component of the boundary  of the cylinders. As
the constructed homeomorphisms coincide on the boundaries, they
set an conjugated homeomorphism of the surface $F.$
\end{proof}

\subsection{Distinguishing graph of the function. }

Let $D = \{F, G_1 (f), $ $...,$ $ G_n(f)\}$ be a diagram of a function
$f$ with isolated critical points on a surface $F.$ Let
cylinder $F_k$ has component of itself boundary on the graphs
$G_i(f), G_{i+1}(f)$. These component of boundary form cycles on
the graphs. The component of boundary and the cycle on
$G_i(f)$ we call by upper and on
$G_{i+1}(f)$ by lower. Thus to local minima and maximas
there correspond cycles consisting of one vertex.
Note that each edge contains in exactly one upper and  one
lower cycle.

If the surface $F$ is oriented then the orientation of surface
induces the orientation on the upper
boundary of each cylinder. This set the orientation of graphs
$G_1, G_2, ..., G_n$. Then upper and lower cycles are oriented.

\begin{defn}
By a distinguishing graph, we call the graphs $G,$ divided on
not intersected subgraphs $G_1, G_2, ..., G_n$ with fixed
upper and lower cycles and the given one-to-one correspondence
between the lower cycles from $G_i$ and upper cycles from
$G_{i+1}, \ i=1, ..., n-1$. Thus each edge contains exactly
in one upper and in one lower cycle, and each isolated vertex
form one upper or lower cycle.
\end{defn}

The isomorphism of the distinguishing graphs is such an
isomorphism of the graphs that  the subgraphs are mapped to
the subgraphs, the upper cycles to the upper ones, lower to lower and
also  the correspondence between the upper and lower cycles
is preserved.

Up to isomorphism, each function $f$ with isolated singular points
on a closed surface sets unique distinguishing graph, which we
call the distinguishing graph of function $f$.

\begin{thm}
 Two functions are topologically conjugate if and only if there
 is an isomorphism of their distinguishing graphs that preserve
 the oder of subgraphs.
\end{thm}

\begin{proof}

 Necessity. If the functions are topologically conjugated,
there is an isomorphism of their diagrams. The restriction of
this isomorphism to the graphs sets isomorphism of the
distinguishing graphs.

Sufficiency. Let distinguishing graphs are isomorphic. We
construct an isomorphism of the diagrams. By the construction,
we have one-to-one correspondence of the cylinders $F_k$ and
the pairs of the upper and lower cycles. Thus the correspondence
between cylinders is given. Let us construct arbitrarily
homeomorphisms between cylinders, so their restrictions on upper
and lower foundation coincide with restrictions of isomorphisms of the
graphs on appropriate cycles. As these homeomorphisms coincide
on the boundaries, they set a homeomorphism of a surface
$F$ and the isomorphism of the diagrams signifies.
\end{proof}

\begin{rem}
If the surface $F$ is oriented and subgraphs have induced
orientations then conjugated homeomorphism preserve the
orientation of the surface only if the isomorphisms of
distinguishing graphs preserve the orientation of cycles.
Correspondent functions are called oriented conjugated.
\end{rem}

It is obvious that functions $f$ and $g$ are topologically
equivalent iff functions $f$ and $g$ or functions $f$ and $-g$
are topologically conjugated.

\begin{cor}
Two functions are topologically equivalent if and only if their
distinguishing graphs are isomorphic.
\end{cor}

\begin{rem}
The distinguishing graphs can be used to classify imbeddings of
graphs into surfaces [3].
\end{rem}

\subsection{Realization}
We discuss problem when the distinguishing graphs
sets function with isolated singular points on a closed surface.

If there is a loop in the graph we can fix a point on it that
divide loop to two edges. We repeat this procedure for each loop.
So  we can suppose that there isn't a loop in the graph.

Two edges $e_1, e_2$ that incident to the same vertex $v$ are
called adjacent for vertex $v,$ if in the distinguishing graph
there is such a cycle, which contains a fragment $(e_1, v, e_2)$
or $(e_2, v, e_1).$

\begin{prop}
A distinguishing graphs sets function with isolated singular
points on a closed oriented surface if and only if for any
two edges that incident to the vertex there is a sequence of
the edges, in which everyone consequent is adjacent in this
vertex with previous.

\end{prop}

\begin{proof}

The necessity is obvious. Let us prove sufficiency. The
conditions of the theorem guarantee that after gluing of
cylinders to the graphs we obtain a closed surface. Let us
construct the function f. For this purpose we set
$f(G_i) =i,$ and for the cylinder $S^1 \times [0,1]$ with
boundary on $G_i$ and $G^{i+1}$: $f(S^1 \times \{t\}) =i+t.$
Thus in vertexes of the graph the cylinders are pasted so, that the function
$f$ coincides with function $i+\text{Re} z^k$ in an appropriate
coordinate system. The constructed function is desired.
\end{proof}

\begin{defn}
A vertex of a distinguishing graph is called planar if the proposition holds for it.
Otherwise it is called conic.
\end{defn}

\section{Examples}

At first we consider function on closed oriented  surfaces with
minimal number of critical points (minimal functions) and
conjugated homeomorphisms that preserve the orientation. On sphere $S^2$
minimal function have two critical points that are minimum and maximum.
It is obvious that all such functions are conjugated.

On other surfaces minimal function have three critical points.
Graphs $G_1, G_3$ are points and graph $G_2$ is a graph with
unique vertex that is a bouquet of $2g+1$ circle, where $g$ is
the genus of surface. Since the cycles on graph $G_1$ and $G_3$
are trivial, the diagram of function consist of graph $G=G_2$
and two cycles on it. Let $v$ be unique vertex of $G$. We orient
the edges of $G$ according to lower cycles and name edges by
$2g+1$ letters: $a, b, c, d, e, ....$  that lower cycles
is given by the word $abcde...$.

Then the upper cycle is an oriented cycle that is described
by word $w$ consisting of the same letters. Thus the function up to
conjugation can be given by this word. Since cyclic perturbation of
the letters in the word lead to the same function, we can
assume that the first letter of the word $w$ is $a$.

Two words correspond to the same function if and only if one can
obtain from other by cyclic renaming of letters ($a\to b, b\to c,
c\to d, ...,\text{last letter} \to a$).

If there is a fragment of two successive letters ($ab,
bc, cd,..$) in the word then the vertex $v$ is conic, because
the end first edge and the ending of the second edge have only
each other as adjective edges. Bellow we consider words without
such fragment.

For torus the words consists of three letters $a, b, c$. It is
obviously there is a unique word $w=acb$ that satisfy condition
above and correspond to graph with planar vertex. Thus up to
oriented conjugation there is a unique  minimal function on the torus.


For surface of genus 2 word consist of letters $a,b,c,d,e$.
There are 4 words:

1) $acbed,$
2) $acebd,$
3) $adbec,$
4) $aedcb.$

The cyclic renaming of letters in each of words 2), 3), 4)
gives the same word and in word 1) all other possible words:
$acedb, adceb, aebdc, $ $aecbd$. All of words 1) - 4) correspond
to the graph with a planar vertex $v$. Thus there are 4 non
oriented conjugated minimal functions on surface of genus 2.

For surface of genus 3 the words consists of 7 letters $a, b,
c, d, e, f, g $. There are 37 words that satisfy the condition
above:

1) $acbedgf,$
2) $acbedgf,$
3) $acbegdf,$
4) $acbfdge,$
5) $acbfegd,$
6) $acbgedf,$
7) $acbgfed,$
8) $acebdgf,$
9) $acebgdf,$
10) $acebgfd,$
11) $acegbdf,$
12) $acfbdge,$
13) $acfbegd,$
14) $acfdbge,$
15) $acfdgeb,$
16) $acfegdb,$
17) $acfebgd,$
18) $acgbfed,$
19) $acgdbfe,$
20) $acgfbed,$
21) $acgfedb,$
22) $adbfcge,$
23) $adbgcfe,$
24) $adbgfec,$
25) $adcgfeb,$
26) $adgcfbe,$
27) $aecbgfd,$
28) $agfedcb,$
29) $aebfcgd,$
30) $afdbgec$

and 31) $acbdgfe,$
32) $acfbged,$
33) $acgebfd,$
34) $acgfdbe,$
35) $adcbgfe,$
36) $adcgfbe,$
37) $aecgfdb.$

The words 1)-30) set planar vertexes, others set conic ones.
For example, if we denote the beginings of the vertexes by
$a_-, b_-, c_-, d_-, e_-,$ $ f_-,$ $ g_-$ and the ends by $a_+, b_+, c_+,
d_+, e_+, f_+, g_+$, then for word 1) we have a chain of
adjacent edges: $a_-c_+b_-e_+d_-g_+ f_- a_+ g_- f_+ e_-
d_+ c_-b_+a_-$. So arbitrary edges (their beginnings and ends)
can be connected by consequence of adjacent. For word 31) we
have three chain: $a_-c_+b_-d_+c_-b_+a_-$, $d_-g_+f_-e_+d_-$ and
$g_-f_+c_-a_+g_-$. The edges from other chain can't be connected
by consequence of adjacent.
Thus there are 30 non oriented conjugated
minimal functions on closed oriented surface of genus 4.

Let us admit that conjugated homeomorphism of surface reverse
the orientation of surface. It correspond to the reverse
orientation of cycles. Then the words $w$ change the order of
letters and change the letters on their reverse. For example,
if word contain 7 letters we change $b-g, c-f,d-e$. So word
$acbegfd$ is changed to $acebdgf$. For surfaces of genus 1 and
2 each word is changed to itself. For surface of genus 3 we
have following 5 pairs of word that are changed one to other:
3) and 8), 2) and 16), 4) and 15), 7) and 21), 17) and 19).
All other word is changed to itself. Denote $k(g)$ the number
of topologically non conjugated minimal functions on the
closed oriented surface of the
genus $g$. So we have: $ k(0)=1, k(1)=1, k(2)=4, k(3)=25.$

Let us consider topological equivalence of minimal functions.
Having the word $w$ for function $f$, we rewrite word $w'$ for
function $-f$. If the word $w=a_1a_2a_3a_4...$, where $a_1, a_2,
a_3,...(a_i \in \{a, b, c, d,... \})$ are letters, then the
word $w$ is obtained from word $abcd...$ by the following
replacement of the letters:
$$a_1 \to a, a_2 \to b, a_3 \to c, a_4 \to d, ...$$

For example in the word $w=acebd$ we replace $a\to a, c\to b,
e\to c, b\to d, d\to e$ in the word $abcde$. Then $w'=adbec$.
So functions given by words $abcde$ and $adbec$ are
topologically equivalent. If $w=acbed$ or $w=aedbc$ then
$w'=w$. Hence there are 3 topologically non equivalent function
on the oriented surface of the genus 2.

For the oriented surface of the genus 3, the following pairs of the
words determinate topologically equivalent functions: $3)-4),
5)-6), 8)-17), 9)-18), 10)-24), 12)-13), 15)-19), 16)-22),
21)-25).$ For all other words $w'=w$. So there are 16
topologically non equivalent functions on the oriented surface
of the genus 3.

Let $f: F \to \Real$ be a minimal function on a non-oriented
surface. As above we denote edges and orient them in such a way
that the lower word have the form $abcd...$. Upper cycle is
non-oriented. The letters that correspond to the edges with
reverse orientation are rewritten with the negative degree. So
functions are given by words consisting of letters with
positive and negative degree.

For $\Real P^2$ there is unique word $ab^{-1}$. For Klein bottle
all words $ab^{-1}c^{-1}, acb^{-1},$ $ac^{-1}b$ are cyclic
equivalent. So there is a unique (up to topological conjugations
and up to topological equivalence) minimal function on $\Real P^2$
and on Klein bottle.

For the non-oriented surface of the genus 3 there are 4 topologically non
conjugated functions that are given by words $1) ab^{-1}c^{-1}d^{-1},$
$2) ab^{-1}dc,$ $ 3) ab^{-1}d^{-1}c,$ $ 4) ab^{-1}dc^{-1}.$ The word
3) and 4) determinate topologically equivalent function. So
there are 3 topologically non equivalent functions on the
non-oriented surface of the genus 3.

The obtained result of the calculations can be summarized in
the table (we write the genus for
the oriented surfaces with "+" and for the non-oriented surfaces
with $"-"$):

\bigskip

\begin{center}

\begin{tabular}{llccccccc}
equivalence $\setminus$ genus &\vline & 0 & 1 & 2 & 3  & $-1$ & $-2$ &
$-3$ \\
\hline
orient. conjugation &\vline & 1 & 1 & 4 & 30 & -  & -  & -  \\
top. conjugation    &\vline & 1 & 1 & 4 & 25 & 1  & 1  & 4  \\
top. equivalence    &\vline & 1 & 1 & 3 & 16 & 1  & 1  & 3
\end{tabular}
\\
the number of non-equivalent minimal functions\\
on the surface of genus $g$
\end{center}

\medskip

If all critical points of function $f$, except of the  minimums and
the maximums, has Poincare index $-1$, then $f$ conjugate with a Morse
function. In [1] the numbers $n(g)$ of topologically non equivalent
Morse functions with one local minimum and one local maximum on
the oriented surface of genus $g \leq 5$ was calculated:
$$ n(1)=1, n(2)=3, n(3)=31, n(4)=778, n(5)=37998.$$

\medskip


\begin{thebibliography}{9}


\bibitem{1}
Kulinich E.V. On topological equivalence of Morse functions on
surfaces // Methods of Func. An. and Topology, No.1, 1998.-
P.22-28.

\bibitem{2}
Kulinich E.V., Prishlyak A.O. On graphs as critical
levels of function on a surface // Some problems of modern math.,
Proc. Inst.Math. Ukrainian NAN, v.25, Kiev. 1998.- P. 87-93.

\bibitem{3}
Prishlyak A.O.  On embedded in surface graphs // Russian
Math. Surveys, v.52, 1997, No.4, 844-845.

\bibitem{4}
Prishlyak A.O.  Vector fields with a given set of singular
points// Ukr.Math.Jorn., v.49, No.10., 1997.- P.1373-1384.

\bibitem{5}
Sharko V.V. On topological equivalence Morse functions on
surfaces // Int. conference at Chelyabinsk State Univ.: Low-
dimensional Topology and Combinatorial Group Theory, 1996.-P.19-
23.

\bibitem{6}
Takens F. The minimal number of critical points of a function
on a compact manifold and Lusternik-Schnirelman category//
Invent. math., No.6, 1968.- P.197-244.
\end{thebibliography}

\end{document}